\documentclass[final,5p,times,twocolumn]{elsarticle}

\usepackage{amssymb}
\usepackage{amsmath}
\usepackage{soul}
\usepackage{subcaption}

\usepackage{doi}
\usepackage{algpseudocode,algorithm,algorithmicx} 

\DeclareMathOperator*{\argmin}{arg\,min}

\usepackage{xcolor}
\newcommand{\od}[1]{#1}
\newcommand{\eb}[1]{#1}

\journal{Journal of Computational Physics}

\begin{document}

\begin{frontmatter}



\title{Nesterov Acceleration for Ensemble Kalman Inversion and Variants}


\author[caltech]{Sydney Vernon} 

\affiliation[caltech]{organization={California Institute of Technology},
            city={Pasadena},
            state={California},
            country={United States}}

\author[caltech,reading,nceo]{Eviatar Bach}
\ead{eviatarbach@protonmail.com}

\affiliation[reading]{organization={University of Reading},
            city={Reading},
            country={United Kingdom}}

\affiliation[nceo]{organization={National Centre for Earth Observation},
            city={Reading},
            country={United Kingdom}}

\author[caltech]{Oliver R. A. Dunbar}
\ead{odunbar@caltech.edu}

\begin{abstract}
Ensemble Kalman inversion (EKI) is a derivative-free, particle-based optimization method for solving inverse problems. It can be shown that EKI approximates a gradient flow, which allows the application of methods for accelerating gradient descent. Here, we show that Nesterov acceleration is effective in speeding up the reduction of the EKI cost function on a variety of inverse problems. We also implement Nesterov acceleration for two EKI variants, unscented Kalman inversion and ensemble transform Kalman inversion. Our specific implementation takes the form of a particle-level nudge that is demonstrably simple to couple in a black-box fashion with any existing EKI variant algorithms, comes with no additional computational expense, and with no additional tuning hyperparameters. This work shows a pathway for future research to translate advances in gradient-based optimization into advances in gradient-free Kalman optimization.
\end{abstract}







\end{frontmatter}


\section{Introduction}

\subsection{Inverse problem}
We assume that we have a model $\mathcal{G}$ with unknown parameters $u^*\in\mathbb{R}^d$, and observations $y\in\mathbb{R}^k$ generated as
\begin{equation}\label{eq:ip}
y = \mathcal{G}(u^*) + \eta,
\end{equation}
where $\eta$ is measurement noise, assumed to be Gaussian $\eta\sim\mathcal{N}(0, \Gamma)$. The inverse problem is then to estimate $u^*$ given $y$. We assume access to the model $\mathcal{G}$. \eb{Note that $y$ can include multiple independent observations of the same quantities.}

We define the cost function \begin{align}
    \mathcal{J}(u)=-\log \mathbb{P}(y|u) = \frac{1}{2}(y - \mathcal{G}(u))^\top \Gamma^{-1} (y - \mathcal{G}(u)), \label{eq:cost}
\end{align}
and note that the maximum likelihood estimate (MLE) of $u^*$ is given by
\begin{equation*}
    u_\text{MLE} = \argmin_u \mathcal{J}(u).
\end{equation*}
\subsection{Ensemble Kalman inversion}
We now turn to the estimation of $u_\text{MLE}$ using ensemble Kalman inversion \cite[EKI,][]{iglesias_ensemble_2013}. \od{While we adopt the notation and terminology of \cite{iglesias_ensemble_2013}, this methodology has a long lineage from randomized likelihood optimizers and ensemble smoothers \cite{EveLee96,EveLee00,CheOliZha09,CheOli12,BocSak14} as used for state and parameter estimation in weather forecasting, and for oil reservoir history matching.} \eb{Moreover, the idea of using an ensemble to approximate a gradient dates back to the introduction of the ensemble Kalman filter \cite{evensen_sequential_1994}.} We note that while variants or extensions \cite{garbuno-inigo_interacting_2020,cleary_calibrate_2021,iglesias_adaptive_2021} of EKI can also be used to approximate the Bayesian inverse problem of estimating a posterior distribution $\mathbb{P}(u|y)$, here we restrict ourselves to finding a point estimate. 

We begin with the continuous-time version of EKI. We first note that
\begin{align}
     -\nabla \mathcal{J}(u) = \frac{d\mathcal{G}}{du}^\top \Gamma^{-1} (y - \mathcal{G}(u)).\label{eq:grad}
\end{align}
Then, we consider having a probability distribution over values of $u$ and define the covariance matrices
\begin{align*}
    C^{uu} &= \mathbb{E}[(u - \mathbb{E}[u])\otimes(u - \mathbb{E}[u])],\\
    C^{u\mathcal{G}} &= \mathbb{E}[(u - \mathbb{E}[u])\otimes(\mathcal{G}(u) - \mathbb{E}[\mathcal{G}(u)])].
\end{align*}
We also note that \eb{under assumptions that the second derivative of $\mathcal{G}$ is small \cite{calvello_ensemble_2024},}
\begin{align}
     C^{u\mathcal{G}} \approx C^{uu}\frac{d\mathcal{G}}{du}\eb{\Bigg|_{\mathbb{E}[u]}^\top}.\label{eq:grad_approx}
\end{align}

We now consider the covariance-preconditioned gradient flow
\begin{align}
    \dot u + C^{uu}\nabla \mathcal{J}(u) = 0.\label{eq:grad_flow}
\end{align}
This differs from a typical gradient flow in the use of the $C^{uu}$ as a time-dependent preconditioner. This is, however, a natural choice, since this type of gradient flow is affine invariant \cite{calvello_ensemble_2024}. Notice that \eqref{eq:grad_flow} is a mean-field or McKean--Vlasov equation, in the sense that the time evolution of $u$ depends on its probability distribution. Given some initial distribution $\mathbb{P}(u(0))$, \eqref{eq:grad_flow} defines a Liouville equation which propagates this distribution forward in time.

Substituting \eqref{eq:grad} and \eqref{eq:grad_approx} into \eqref{eq:grad_flow}, we obtain
\begin{align}
    \dot u = C^{u\mathcal{G}} \Gamma^{-1} (y - \mathcal{G}(u)).\label{eq:mf_eki}
\end{align}
This is the mean-field EKI \cite{calvello_ensemble_2024}. Importantly, the derivative of the model $\frac{d\mathcal{G}}{du}$ is replaced with the derivative-free covariance approximation. We now approximate this further, and replace $C^{u\mathcal{G}}$ by its Monte Carlo approximation from a finite ensemble of particles $\{u^{(n)}\}_{n=1}^N$ and discretize in time with an explicit method using a time-step of $\Delta t$, obtaining the usual form of EKI:
\begin{subequations}
\begin{align}
    C^{u\mathcal{G}}_j &= \frac{1}{N}\sum_{n=1}^N \left[ (u^{(n)}_j - \overline{u_j})\otimes(\mathcal{G}(u^{(n)}_j) - \overline{ \mathcal{G}(u)}_j) \right],\\
    C^{\mathcal{GG}}_j &= \frac{1}{N}\sum_{n=1}^N \left[(\mathcal{G}(u^{(n)}_j) - \overline{ \mathcal{G}(u)}_j)\otimes(\mathcal{G}(u^{(n)}_j) - \overline{ \mathcal{G}(u)}_j) \right],\\
    \overline{u}_j &= \frac{1}{N}\sum_{n=1}^N u^{(n)}_j,\qquad \overline{\mathcal{G}(u)}_j = \frac{1}{N} \sum_{n=1}^N\mathcal{G}(u^{(n)}_j),\\
    u_{j+1}^{(n)} &= u_j^{(n)} + \Delta t C^{u\mathcal{G}}_j(\Gamma + \Delta t C_j^{\mathcal{GG}})^{-1}(y - \mathcal{G}(u_j^{(n)})). \label{eq:eki-alg}
\end{align}
\end{subequations}
Note that this discretization is not quite a forward Euler method, as in \eqref{eq:eki-alg} we have replaced $\Gamma^{-1}$ by $(\Gamma + \Delta t C_j^{\mathcal{GG}})^{-1}$. This discretization, however, is still consistent with the limit of \eqref{eq:mf_eki} as $\Delta t\to 0$, and also maintains consistency with a derivation of EKI based on the discrete-time ensemble Kalman filter as well as having favorable convergence properties; see \cite{blomker_strongly_2018}.

Sometimes, the initial ensemble $\{u^{(n)}_0\}_{n=1}^N$ is drawn from a prior distribution over $u$. Since we are here interested in the MLE, and due to the invariant subspace property of EKI \cite{iglesias_ensemble_2013,tong_localized_2023}, the choice of initial ensemble can be considered a regularization for the maximum likelihood problem; that is, EKI approximates the MLE but with regularization coming from the fact that the solution must be in the span of the initial ensemble. The recursion \eqref{eq:eki-alg} is repeated until some $j = J$, and the ensemble mean $\overline{u}_J$ taken as an approximation of $u_\text{MLE}$. The convergence of the particle approximation to the mean-field limit is proved in \cite{ding_ensemble_2021}.

\eb{We also note that EKI is often derived from the ensemble Kalman filter or from connection with the Bayesian solution of a linear inverse problem, wherein observations are perturbed for each ensemble member \cite{iglesias_ensemble_2013}. However, in this paper, since we use the formulation of EKI based on MLE, perturbed observations do not arise.}

Several variants of EKI exist. Unscented Kalman inversion \cite[UKI,][]{huang_iterated_2022} chooses the ensemble using quadrature points at every iteration and has favorable performance for some problems. Ensemble transform Kalman inversion (ETKI), based on the ensemble transform Kalman filter \cite{tippett_ensemble_2003}, avoids building covariance matrices in the parameter and observation spaces, and has favorable computational complexity if $N\ll d,k$. 

\subsection{Nesterov acceleration}

Nesterov acceleration, sometimes referred to as momentum, is a method to accelerate gradient descent \cite{nesterov_method_1983,nesterov_introductory_2004}. This was extended to continuous-time gradient flows in \cite{su_differential_2016}. Namely, continuous-time Nesterov acceleration transforms a gradient flow
\begin{equation*}
    \dot u + \nabla\mathcal{J}(u) = 0
\end{equation*}
to
\begin{align}
\ddot u + \lambda(t) \dot u + \nabla \mathcal{J}(u) = 0.\label{eq:nesterov}
\end{align}
Here, $\lambda(t)$ is the momentum coefficient, which is taken in \cite{su_differential_2016} to be $\lambda(t) = \frac{3}{t}$. In this paper, we will apply this Nesterov acceleration approach to EKI, as well as its variants UKI and ETKI.

\subsection{Motivation and related work}

Each iteration of \eqref{eq:eki-alg} requires evaluating the forward model $\mathcal{G}(\cdot)$ $N$ times. This is often computationally expensive, as in the case of climate models \cite{cleary_calibrate_2021}. Thus, we would like $J$ to be as small as possible, which motivates the need for acceleration.

Nesterov acceleration for EKI was first suggested by \cite{kovachki_ensemble_2019}; however, detailed numerical tests were not conducted in that work, nor was Nesterov applied to EKI variants. Another method for accelerating convergence of EKI was considered in \cite{chada_convergence_2022}. \od{A related implementation was also considered for the ensemble Kalman updates with a different objective \cite{nilsen_accelerated_2024}. There momentum was applied to the ensemble mean and covariance, and with constant momentum parameter. By contrast, we will showcase a simpler particle-wise update, and one with a momentum parameter that changes with iteration in a manner consistent with the gradient-based theory.} Momentum has also been considered for particle-based sampling methods in \cite{liu_second_2022}.

\section{Nesterov acceleration for EKI}

We now add Nesterov acceleration to EKI. We begin by considering a covariance-preconditioned Nesterov gradient flow obtained by formally replacing the gradient in \eqref{eq:nesterov} with the covariance-preconditioned gradient as in \eqref{eq:grad_flow},
\begin{equation*}
\ddot u + \lambda(t)\dot u + C^{uu}\nabla \mathcal{J}(u) = 0.
\end{equation*}
We note that although improved rates of convergence have previously been proved for continuous-time Nesterov, it remains to be proved that the covariance-preconditioned Nesterov will also have this improvement; however, our numerical results suggest this to be the case. 

Now, employing \eqref{eq:grad} and \eqref{eq:grad_approx}, we obtain the mean-field Nesterov-accelerated EKI
\begin{align}
\ddot u + \lambda(t)\dot u = C^{u\mathcal{G}} \Gamma^{-1} (y - \mathcal{G}(u)).\label{eq:nesterov_eki}
\end{align}

Due to the non-uniqueness of discretization of \eqref{eq:nesterov_eki}, it is difficult to proceed directly with the derivation from these mean-field equations. Instead, we next write our proposed discrete algorithm, and then motivate its consistency with \eqref{eq:nesterov_eki} by taking small-timestep and large-particle limits.


\begin{algorithm}
\caption{Nesterov-accelerated ensemble Kalman inversion}
\label{alg:nest-eki}
\begin{algorithmic}[1]
\Require $\{u_0^{(n)}\}_{n=1}^N$, $J \in \mathbb{N}$, $y$, $\Gamma$, $\mathcal{G}$
\State Compute $C_1^{u\mathcal{G}},C_1^{\mathcal{G}\mathcal{G}}$
\State $u_1^{(n)} \gets u_0^{(n)} + \Delta t C^{u\mathcal{G}}_1(\Gamma + \Delta t C_1^{\mathcal{GG}})^{-1}(y - \mathcal{G}(u_1^{(n)})),\, \forall n=1,\dots,N$ 
\For{$j=1,\dots,J$}
\State $\lambda_j \gets \frac{j-1}{j+2},$
\State $v^{(n)}_j \gets u^{(n)}_j + \lambda_j (u^{(n)}_j -u^{(n)}_{j-1}), \forall n=1,\dots,N$
\State Compute $C_j^{v\mathcal{G}},C_j^{\mathcal{G}\mathcal{G}}$
\State $u^{(n)}_{j+1} \gets v^{(n)}_j + \Delta tC_{j}^{v\mathcal{G}} (\Gamma + \Delta t C^{\mathcal{G}\mathcal{G}}_j)^{-1} (y-\mathcal{G}(v^{(n)}_j)),\forall n=1,\dots,N$
\EndFor
\State \Return $\{u^{(n)}_{J+1}\}_{n=1}^N$
\end{algorithmic}
\end{algorithm}
We write out the Nesterov-accelerated EKI as Algorithm~\ref{alg:nest-eki}. One can view this implementation of Nesterov acceleration as performing a particle-wise ``nudge'' based on each particle's history, and then carrying out the regular EKI update \eqref{eq:eki-alg} with the nudged ensemble. 

In Algorithm~\ref{alg:nest-eki} we define $\lambda_j$ in line with the Nesterov acceleration literature; however, this choice is not unique. We explore other choices and their effect on the convergence in \ref{sec:lambdas}. 

With any choice of $\lambda_j$, the accelerated algorithm retains the subspace property of EKI, that particles will always remain in the span of the initial ensemble \cite{iglesias_ensemble_2013}. This can be seen inductively, from the steps of Algorithm~\ref{alg:nest-eki}. Firstly, $\{u^{(\cdot)}_0\}$ and $\{u^{(\cdot)}_1\}$ are by definition in this span; then, if $\{u^{(\cdot)}_{k-1}\}$, $\{u^{(\cdot)}_k\}$ are in the span, $\{v^{(\cdot)}_k\}$ is too, as each member update is a linear combination of elements in the span. Then by the original argument \cite{iglesias_ensemble_2013}, $\{u^{(\cdot)}_{k+1}\}$ is also in the span, as a classical EKI update of $\{v^{(\cdot)}_k\}$.

In \ref{sec:cts-limit} we use a formal argument to find the continuous-time limit of Algorithm~\ref{alg:nest-eki}, but highlight a need for stronger rigorous arguments to show convergence to \eqref{eq:nesterov_eki} in the limit of infinite ensemble size.

Applying the same idea to the UKI and ETKI algorithms, we obtain Nesterov-accelerated versions of these two algorithms in \ref{sec:uki_etki}.

\section{Numerical experiments}\label{sec:exp}

We conduct numerical experiments on the impact of Nesterov acceleration on three different inverse problems. Data is generated from equation \eqref{eq:ip}, so we must define the forward map $\mathcal{G}$, noise covariance $\Gamma$, and select a ground-truth parameter $u^*$. We must additionally pick a probability distribution on $u_0$, from which the initial ensemble is drawn. We present results of experiments over fifty trials, each with different random draws of the initial ensemble and different noise realizations $\eta$ in \eqref{eq:ip}. When reporting the cost function values in the experiments that follow, we use the ensemble mean over the forward map evaluations $\overline{\mathcal{G}(u)}_j$ in expression \eqref{eq:cost}.

\subsection{Example inverse problems}\label{sec:ex-ip}


In the exponential sine problem (Exp Sin), we consider:
\begin{align*}
    u &= [u_1, u_2]^\top,\\
    f(t,u) &= \exp\left(u_1 \sin(t) + u_2\right),\\
    \mathcal G(u) &= \begin{bmatrix}
        \mathrm{mean}_{t \in [0,2\pi]}\left(f(t,u)\right) \\
        \max_{t \in [0,2\pi]}(f(t,u)) - \min_{t \in [0,2\pi]}(f(t,u))
        \end{bmatrix}.
\end{align*}
That is, we look to estimate the amplitude and vertical shift parameters $u = \left[ u_1, u_2 \right]^\top$ of function $f$, using the observations given by $\mathcal G$.  We define the true parameters to be $u^* = \left[1, 0.8\right]^\top$. We take a Lognormal$ (-1.38,0.06)$ initial distribution for $u_1$ and $\mathcal{N}(0,0.5)$ for $u_2$, and assume these are independent.

In the Lorenz '96 problem (Lorenz96), we solve for the initial conditions of the $D$-dimensional Lorenz '96 dynamical system \cite{lorenz_predictability:_1996} given a state some time after that is integrated from those initial conditions, following a problem used by \cite{tong_localized_2023}. 
Here, we have
\begin{align*}
    u &= x(t_0) \in \mathbb R^D, \qquad   \mathcal G(u) = x(t_K) \in \mathbb R ^ D, \\
    \frac{dx_{k}}{dt} &= -x_{k} - x_{{k-1}} (x_{{k-2}} - x_{{k+1}}) + F,
\end{align*}
where the indices are cyclical and $x(t) = [x_1(t), \ldots, x_D(t)]^\top$. 
We use the fourth-order Runge--Kutta method with timestep $0.05$, and $t_K = 0.4$. The dimension is set to $D = 20$ and the forcing parameter to $F = 8$, for which the Lorenz '96 model is chaotic. 
We draw the truth $u^*$ from the $D$-dimensional standard normal distribution.
Taking $u^*$ to be the true state at time $t_0$, we
run through the dynamics until time $t_k = t_0 + 1000$.
The initial ensemble members for $u$ are each drawn from the $D$-dimensional standard normal distribution as well.

Finally, the Darcy problem (Darcy) describes subsurface flow over the domain $\mathcal D = [0,1]^2$. We observe pointwise measurements from a pressure field $p(x)$. We seek to optimize parameter $\kappa$ to define a discretized permeability field $a(x, \kappa)$. Due to the high dimensionality of the discretized field when spatial resolution is high, we choose to represent it efficiently with a Gaussian random field (GRF) of a given covariance structure (e.g., see the literature review of \cite{EisLatUll19}). We choose a Matérn covariance of smoothness $1$ and correlation length $0.25$. This covariance is expanded as a Karhunen--Loève eigenfunction expansion, and truncated to $d$ terms (degrees of freedom); the learned $u$ are then coefficients of the truncated expansion.
Our forward map then becomes:
\[
    u \in \mathbb R^{d}, \quad \mathcal G(u) = \mathrm{vec}\left(  
    \begin{bmatrix}
    p(x_{1,1}) & \cdots & p(x_{1,N}) \\
    \vdots & \ddots & \vdots \\
    p(x_{N,1}) & \cdots & p(x_{N,N})
    \end{bmatrix}\right),
\]
\noindent where $p(x)$ is obtained by solving Darcy's law with a finite-difference scheme. Darcy's law is given by:
\begin{align*}
    -\nabla \left( a(x,\kappa(u)) p(x) \right) = f(x),
\end{align*}
where $f(x)$ is a given velocity field (flow) and $p(x)$ is set to $0$ at the boundary. We set $d = 50$ and $N$ = 80, and $u^* =[ -1.5,\dots,-1.5]^\top$. The distribution of the initial ensemble for $u$ is taken to be independent $\mathcal{N}(0,1)$ for each degree of freedom.


\subsection{Implementation}

We conduct experiments using the open-source Julia package \texttt{EnsembleKalmanProcesses.jl} \cite{dunbar_ensemblekalmanprocessesjl_2022}. This package has implementations of EKI, UKI, and ETKI algorithms, and we have added Nesterov acceleration implementations for these variants as well as the different choices of $\lambda_j$.

\subsection{Results}

\begin{figure*}
    \centering

    \vspace{-20pt}
    
    \includegraphics[width=0.87\linewidth]{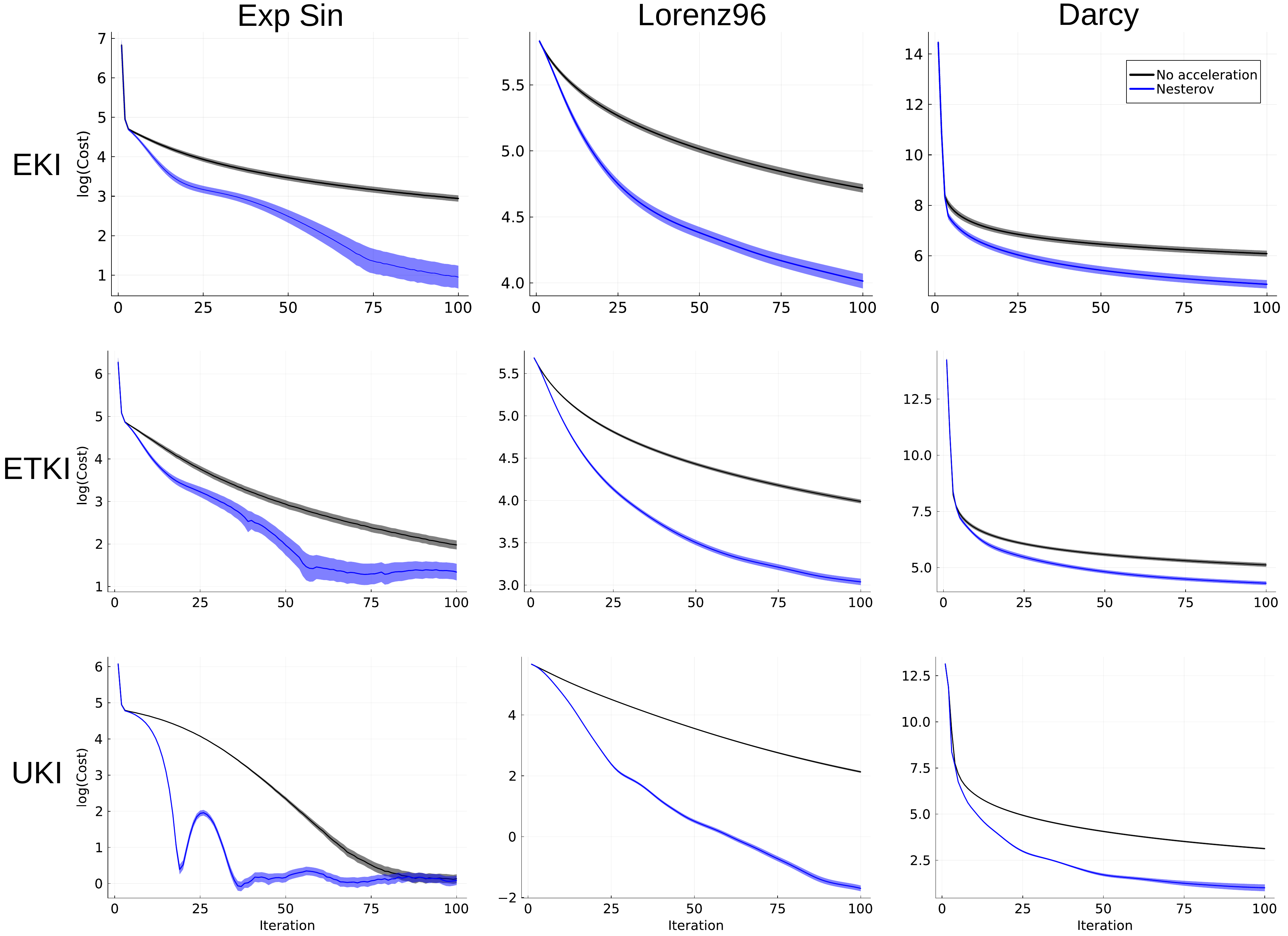}

    \vspace{-10pt} 
    
    \caption{\eb{Comparison of EKI, ETKI, and UKI convergence with and without acceleration, over the three inverse problems of Section~\ref{sec:ex-ip}, as measured by the logarithm of the cost function \eqref{eq:cost}. 
    The results were averaged over 50 trials for exponential sine and Lorenz '96, and 10 trials for Darcy. Ribbons denote one standard error from the mean.}}
    \label{fig:nesterov_results}
\end{figure*}


\od{Figure~\ref{fig:nesterov_results} shows the results of applying Nesterov acceleration to the three algorithms EKI, ETKI, and UKI on the three inverse problems. For the first and second row, EKI and ETKI, the ensemble size is taken to be 10, 20, and 52, respectively, for the three inverse problems, weakly related to the increasing problem dimension (\ref{sec:additional_experiments} shows that acceleration occurs regardless of ensemble size). The results are summarized over random trials where the data realization $y$ and the initial ensemble $\{u^{n}_0\}_{n=1}^N$ are both redrawn. In all three cases, and for all three algorithms, Nesterov acceleration results in faster reduction of the cost function, as well as appearing to converge to better optima. This suggests that Nesterov acceleration can successfully be applied to EKI variants, such as those described in \cite{huang_iterated_2022}.}

   


For additional verification of the robustness of the method, we conducted three other experiments; the results are shown in \ref{sec:additional_experiments}. Namely, we test the performance of Nesterov-accelerated EKI with a variety of ensemble sizes, with different timesteps $\Delta t$, and with different expressions for the momentum coefficient. In all cases, Nesterov-accelerated EKI exhibits improvement compared to non-accelerated EKI.

In all the experiments, Nesterov-accelerated EKI performed at least as well as non-accelerated EKI, and often significantly better, across problems with varying dimensionality of the input space and nonlinearity of the forward model. Despite some variance over different random trials, the improvements are robust also when looking at individual trials. In all the cases tested, the improvement in the decrease of the cost function was observed over the first few iterations in addition to performing better in the long run, indicating that Nesterov acceleration can be useful even in settings where only a small number of iterations $J$ can be afforded.

In some of the cases, such as the UKI experiment \od{(bottom row of Figure~\ref{fig:nesterov_results})}, temporary ``bumps'' where the cost function increases with Nesterov-accelerated EKI can be observed. This may be because Nesterov-accelerated gradient descent, unlike gradient descent, is not strictly a descent method, in the sense that the objective function \eqref{eq:cost} is not guaranteed to be monotonically non-increasing with iterations. This behavior is also observed when accelerating gradient descent \cite{su_differential_2016}.

\section{Conclusions}

We have demonstrated that Nesterov acceleration improves performance of EKI on a variety of inverse problems. The form of update we used here is written as a parameter-wise nudging that is simple to implement, non-intrusive to the EKI update, and instantly adaptable to other algorithm variants. There is no additional computational cost, nor are there additional algorithm hyperparameters to tune. We demonstrated this with algorithms for accelerated UKI and ETKI, and also observed improved convergence.

\section{Future work}
On a theoretical level, it remains to be proven that the Nesterov acceleration in fact improves the convergence rate of covariance-preconditioned gradient flows, as has previously been done for regular gradient flows. Adaptive timestepping through modification of $\Delta t$ has been explored for EKI \cite{iglesias_adaptive_2021}, and it remains for future work to apply acceleration in such settings. The motivation of \ref{sec:cts-limit} also still requires rigorous justification to show whether \eqref{eq:nesterov_eki} is the mean-field limit of \eqref{eq:ens-flow}.

Another avenue of future work is the extension of the current work to the ensemble Kalman sampler \cite[EKS,][]{garbuno-inigo_interacting_2020}. Although there has been previous work on acceleration of continuous-time EKS \cite{wang_accelerated_2021,liu_second_2022}, the discrete-time versions of the algorithms are complex and depend in detail on the form of the EKS update. Our approach may offer simple implementations of such updates, and \od{more broadly,  application to other empirically approximated non-Bayesian objective functions.}

This work demonstrates that improvements of gradient descent algorithms can be transferred to EKI and its variants, by using the EKI's continuous-time gradient flow limit as a bridge. This paves the way for other gradient-based optimizer advancements popular with the machine learning community  to be transferred into the EKI framework. \od{Examples of recently developed methods that could be implemented for EKI acceleration are the Information Theoretic Exact Method (ITEM) \cite{DroTay23} or its cousin, triple momentum \cite{VanSvoy_2017,AspSciTay21} (corresponding to ITEM with constant momentum parameters); these methods have provable optimality. Additionally, there is the popular family of adaptive methods (namely, Adam, Nadam, Adan, Adamax, and Adagrad \cite{kingma_adam_2017,dozat_incorporating_2016,xie_adan_2024}) that lack in theoretical guarantees and continuum limits, but are used commonly in the machine learning community.}

\section*{Acknowledgments}

SV was supported through a Summer Undergraduate Research Fellowship (SURF) at Caltech. EB was supported by the Foster and Coco Stanback Postdoctoral Fellowship at Caltech. EB and ORAD were supported by the Schmidt Foundation and by the Office of Naval Research (Grant No. N00014-23-1-2654). We thank Andrew Stuart for helpful suggestions.
\appendix
\section{Continuous-time limit for Nesterov-accelerated EKI}\label{sec:cts-limit}
This analysis closely follows the limit of \cite{su_differential_2016}.
First, for brevity define 
\begin{align*}
\mathcal{F}_n(\{w_j^{(\cdot)}\}) &= C_{j}^{u\mathcal{G}}(\{w_j^{(\cdot)}\}) (\Gamma + \Delta t C^{\mathcal{G}\mathcal{G}}(\{w_j^{(\cdot)}\}))^{-1} (y-\mathcal{G}(w^{(n)}_j)),\\
\overline{\mathcal{F}}_n(\{W^{(\cdot)}(t)\}) &= C^{u\mathcal{G}}(\{W^{(\cdot)}(t)\}) \Gamma^{-1}(y-\mathcal{G}(W^{(n)}(t))).
\end{align*}
These are the discrete- and continuous-time expressions for the ensemble update, respectively. We use the notation $\{x^{(\cdot)}\}$ to denote the full ensemble $\{x^{(n)}\}_{n=1}^N$. We also assume that $\mathcal{G}$ is a continuous function of its argument. $C_{j}^{u\mathcal{G}}(\{w_j^{(\cdot)}\})$ and $C^{\mathcal{G}\mathcal{G}}(\{w_j^{(\cdot)}\})$ indicate the respective empirical covariances computed using the ensemble. We can express the discrete Nesterov update of Algorithm~\ref{alg:nest-eki} as follows:
\begin{align}
    v^{(n)}_j &= u^{(n)}_j + \lambda_j  (u^{(n)}_j - u^{(n)}_{j-1}),\nonumber\\
    u^{(n)}_{j+1} &= v^{(n)}_j - \Delta t \mathcal{F}_n(\{v_j^{(\cdot)}\}).\label{eq:disc_nest}
\end{align}
Rewriting and using $\lambda_j = \frac{j-1}{j+2}$, this becomes
\begin{equation}
\label{eq:dt_disc_nest}
    \frac{u^{(n)}_{j+1} - u^{(n)}_j}{\!\sqrt{\Delta t}} = \frac{j-1}{j+2}\left(\frac{u^{(n)}_j-u^{(n)}_{j-1}}{\!\sqrt{\Delta t}}\right) - \!\sqrt{\Delta t} \mathcal{F}_n(\{v^{(\cdot)}_j\}).
\end{equation}

We associate a continuous-time trajectory $U^{(n)}(t)$ with the discrete trajectories $u^{(n)}_j$ via the ansatz $u^{(n)}_j = U^{(n)}(t) = U^{(n)}((j-2)\sqrt{\Delta t})$, that is, set $j = t/\sqrt{\Delta t} - 2$. We assume that $U^{(n)}(t)$ is twice-differentiable with respect to $t$. Then, one can Taylor expand $U(t + \sqrt{\Delta t})$ and $U(t - \sqrt{\Delta t})$ about $t$ to derive the expressions:
\begin{align*}
   \frac{u^{(n)}_{j+1} - u^{(n)}_{j}}{\sqrt{\Delta t}} &= \dot U^{(n)}(t) + \frac{1}{2}\sqrt{\Delta t} \ddot U^{(n)}(t) + o(\!\sqrt{\Delta t}),\\
   \frac{u^{(n)}_{j} - u^{(n)}_{j-1}}{\sqrt{\Delta t}} &= \dot U^{(n)}(t) - \frac{1}{2}\sqrt{\Delta t} \ddot U^{(n)}(t) + o(\!\sqrt{\Delta t}).
\end{align*}
Taylor expansion of $U(t - \sqrt{\Delta t})$ in $v_j^{(n)}$ about $t$ gives:
\begin{align}
v^{(n)}_j &= U^{(n)}(t) + \frac{j-1}{j+2}\left(U^{(n)}(t) - U^{(n)}(t-\!\sqrt{\Delta t})\right), \nonumber\\
& = U^{(n)}(t) + \!\sqrt{\Delta t}\frac{j-1}{j+2}\dot U(t)+o(\!\sqrt{\Delta t}).\label{eq:v_about_t}
\end{align} 

We would now like to show that we can write
\begin{equation}\label{eq:Fv_FU}
    \mathcal{F}_n(\{v_j^{(\cdot)}\}) = \mathcal{F}_n(\{U^{(\cdot)}(t)\})+ o(1).
\end{equation}
To do so, we first note that for matrices $\eb{A} = \Gamma$ and $\Delta t \eb{B} = -\Delta t C^{\mathcal{G}\mathcal{G}}(U^{(\cdot)}(t))$, and taking $\Delta t$ small enough so that $\Delta t\|\eb{A}^{-1}\eb{B}\|<1$, where $\|\cdot\|$ is the operator norm, it follows from Neumann series expansion that 
\[
(\eb{A}-\Delta t \eb{B})^{-1} = \eb{A}^{-1}+\sum_{k=0}^\infty (\Delta t)^{k}(\eb{A}^{-1}\eb{B})^k\eb{A}^{-1} = o(\Delta t).
\]
We then substitute \eqref{eq:v_about_t} into $\mathcal{F}_n(\{v_j^{(\cdot)}\})$ and Taylor expand to first order:
\begin{align}
    &\left(C_{j}^{u\mathcal{G}}(\{U^{(\cdot)}(t)\}) + \sqrt{\Delta t}c_1 + o(\sqrt{\Delta t})\right)\label{eq:f_v}\\
    &\qquad\times \left(\Gamma + \Delta t C^{\mathcal{G}\mathcal{G}}(\{U^{(\cdot)}(t)\}) + (\Delta t)^{3/2} c_2 + o((\Delta t)^{3/2})\right)^{-1}\nonumber\\
    &\qquad\times\left(y-\mathcal{G}(U^{(n)}(t)) + \sqrt{\Delta t}c_3 + o(\sqrt{\Delta t}) \right),\nonumber
\end{align}
where $c_1$, $c_2$, and $c_3$ are not functions of $\Delta t$. Then, defining $\eb{D} = \Gamma + \Delta t C^{\mathcal{G}\mathcal{G}}(\{U^{(\cdot)}(t)\})$, we have for sufficiently small $\Delta t$ the Neumann series expansion 
\begin{align*}
    &(\eb{D} + (\Delta t)^{3/2} c_2 + o((\Delta t)^{3/2}))^{-1}\\
    &\qquad= \left(I + \eb{D}^{-1}(\Delta t)^{3/2} c_2 + o((\Delta t)^{3/2})\right)^{-1}\eb{D}^{-1},\\
    &\qquad= (I - \eb{D}^{-1}(\Delta t)^{3/2}c_2 + o((\Delta t)^{3/2}))\eb{D}^{-1},
\end{align*}
where we have used the fact that $\eb{D}^{-1}$ is of order $\Delta t$ for sufficiently small $\Delta t$. Substituting back into \eqref{eq:f_v} and expanding, we find the desired result \eqref{eq:Fv_FU}.


Using \eb{again the fact that} $\eb{D}^{-1}$ is $o(\Delta t)$ \eb{for $\Delta t$ sufficiently small}, we deduce
\[
\!\sqrt{\Delta t}\mathcal{F}_n(\{U^{(\cdot)}(t)\}) = \!\sqrt{\Delta t}\,\overline{\mathcal{F}}_n(\{U^{(\cdot)}(t)\})+ o(\!\sqrt{\Delta t}).
\]
Substitution into \eqref{eq:dt_disc_nest} results in
\begin{align*}
    \dot U(t) + \frac{1}{2}\sqrt{\Delta t} \ddot U(t) + o(\!\sqrt{\Delta t}) &= \left(1-\frac{3\sqrt{\Delta t}}{t}\right) \left(\dot U(t) - \frac{1}{2}\sqrt{\Delta t} \ddot U(t) \right)\\
    &\qquad+ \sqrt{\Delta t}\,\overline{\mathcal{F}}_n(\{U^{(\cdot)}(t)\})+ o(\!\sqrt{\Delta t}).
\end{align*}
By asymptotic matching on the order $\sqrt{\Delta t}$, one derives an informal relationship:
\begin{equation}\label{eq:ens-flow}
\ddot U(t) +\frac{3}{t} \dot U(t) = -\overline{\mathcal{F}}_n(\{U^{(\cdot)}(t)\}).
\end{equation}
From here, one would like to take the large-particle limit $N\to \infty$ and show that the mean-field equation is of the form \eqref{eq:nesterov_eki}. As seen in \cite{ding_ensemble_2021}, where this is carried out for the standard EKI algorithm, this ends up being a technical argument for even linear or mildly nonlinear $\mathcal{G}$, and theory involves jointly taking limits $\Delta t \to 0$ and while fixing $N\Delta t = 1$. Here we have taken $\Delta t \to 0$ before $N\to\infty$ as a motivation, but these limits may not commute.

\section{Choice of $\lambda_j$}\label{sec:lambdas}
Different choices of $\lambda_j$ (respectively $\lambda(t)$) are made in the literature. In this work we considered three options:
\begin{itemize}
    \item \textbf{Original} \cite{nesterov_introductory_2004}: $\lambda_j = \frac{j-1}{j+2}$.
    \item \textbf{Recursive} \cite{su_differential_2016}: $\lambda_j = \theta_j(\theta_{j-1}^{-1}-1)$, 
    where $\theta_0=1$,\\  $\theta_{j+1} = \frac{1}{2}\left(\!\sqrt{\theta_j^4 + 4\theta_j^2} - \theta_j^2\right)$.
    \item \textbf{Constant}: $\lambda_j\equiv c$, with $c$ a tuning parameter.
\end{itemize}
Both the Original and Recursive definitions asymptotically behave as $1-\frac{3}{j}+ O\left(\frac{1}{j^2}\right)$, a property needed for acceleration of \eqref{eq:nesterov}; see \cite{su_differential_2016}. The Constant definition, on the other hand, has been shown not to improve asymptotic convergence rates for \eqref{eq:nesterov}; see \cite{KovStu21}. It is, however, a common practice in the literature. In the numerical experiments that follow in Section \ref{sec:exp}, we used the Recursive definition unless otherwise stated. 

The difference in coefficients, and their effect on performance, is illustrated in a numerical example below.

\begin{figure}[h]
    \centering
 \includegraphics[width=0.99\linewidth]{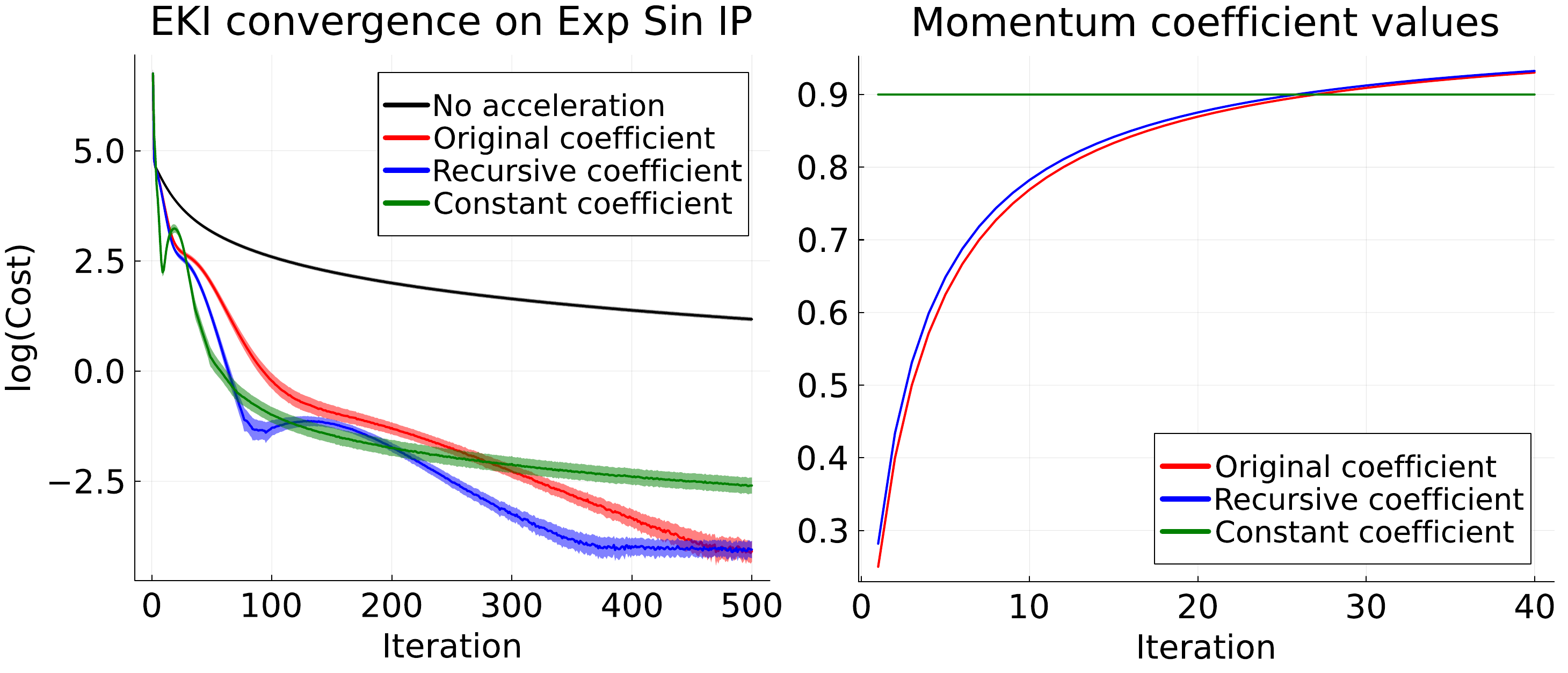}\
    
    \caption{The effect of different choices of the $\lambda_j$ coefficient for EKI acceleration on the Exp Sin inverse problem of Section\ref{sec:ex-ip}. Left: convergence as measured with the logarithm of the cost function \eqref{eq:cost} Ribbons denote one standard error from the mean. Right: values of the $\lambda_j$ over the first 40 iterations.}
    \label{fig:nesterov_other_coefficients}
\end{figure}
Figure~\ref{fig:nesterov_other_coefficients} shows the results of applying the two different expressions for the momentum coefficient on the exponential sine problem, as well as a fixed momentum coefficient $\lambda_j = 0.9$ that was tuned for good performance on this problem. We observe both the Original and Recursive give consistently better performance than without acceleration, and moreover appear to attain a higher rate of convergence. In this experiment we see that the Recursive slightly outperforms the Original, though this is likely experiment dependent. The Constant coefficient actually performs worse than without acceleration for the first 10 iterations, then achieves fast convergence for 100 or so iterations, before its convergence stagnates and eventually the dynamic coefficients overtake it. 

\section{Ensemble transform Kalman inversion and unscented Kalman inversion}\label{sec:uki_etki}

We present the algorithms for accelerated ensemble transform Kalman inversion (ETKI) in Algorithm \ref{alg:nest-etki} and accelerated unscented Kalman inversion (UKI) in Algorithm~\ref{alg:nest-uki}. In both cases, the only difference in the accelerated versions of the algorithms is the introduction of the ``nudged'' ensembles $\{v^{(n)}_j\}_{n,j}$ over which we evaluate the model $\mathcal{G}$. We still write out the full algorithms in easily implementable forms, as they often look significantly different from other forms written in the existing literature (e.g., \cite{huang_iterated_2022}). These forms also highlights the simplicity of the proposed accelerations.

We do not provide derivations of the new methods here but offer some motivation and interpretation. To motivate ETKI, we recall that one can derive the EKI from the stochastic ensemble Kalman filter (EnKF) \cite{iglesias_ensemble_2013}. ETKI is found by applying a similar derivation to the class of square-root filters known as ensemble transform Kalman filter (ETKF); the power of this update differs from EKI, in that one computes all matrices in the ensemble space, allowing for preferential scaling to EKI in the output dimensions, assuming $\Gamma$ is easily invertible \cite{tippett_ensemble_2003}. It also does not accumulate errors due to stochastic perturbations.

UKI is obtained, in its basic form, by replacing EKI's Monte Carlo approximation of a Gaussian with a known quadrature rule. It also contains additional regularization, typically by the prior. In lower dimensions this can lead to significant efficiency increases, and very fast convergence on some problems. There are a number of free parameters and matrices in this algorithm, and for our experiments we select these based on the work of \cite{huang_iterated_2022}: we take the prior to be $N(m,C)$ and take $r=m$, $\Sigma_\nu=2\Gamma$, and $\Sigma_\omega = (2-\alpha^2)C$. We therefore need only to choose $\alpha \in (0,1]$ in experiments. UKI picks the initial ensemble deterministically; therefore, the different trials in the numerical experiments vary only due to using different data realizations.

\begin{algorithm}
\caption{Nesterov-accelerated ensemble transform Kalman inversion}
\label{alg:nest-etki}
\begin{algorithmic}[2]
\Function{$\mathrm{compute\_increment}$}{$u,\overline{u},\{\mathcal{G}(u^{(n)})\}_n,\Gamma$}
\State Compute $\overline{\mathcal{G}}(u) \gets \frac{1}{N}\sum_n\mathcal{G}(u^{(n)})$
\State $\mathrm{col}_n(\Delta u) \gets \frac{1}{\sqrt{N-1}}(u^{(n)}-\overline{u})$
\State $\mathrm{col}_n(\Delta \mathcal{G}(u)) \gets \frac{1}{\sqrt{N-1}}(\mathcal{G}(u^{(n)})-\overline{\mathcal{G}}(u))$
\State $\Omega(u) \gets (I + \Delta \mathcal{G}(u)^{\top} \Gamma^{-1} \Delta \mathcal{G}(u))^{-1}$
\State $w(u) \gets \Omega \Delta \mathcal{G}(u)^{\top}\Gamma^{-1} (y - \overline{\mathcal{G}}(u))$
\State $\delta u \gets \overline{u} + \Delta u (w(u) + \sqrt{N-1}\sqrt{\Omega(u)})$ 

\State \Return $\delta u$
\EndFunction
\State
\Require $u_0 = \{u_0^{(n)}\}_{n=1}^N$, $J \in \mathbb{N}$, $y$, $\Gamma$, $\mathcal{G}$
\State Compute $\overline{u}_0 \gets \frac{1}{N}\sum_n u_0^{(n)}$
\State $\delta u_0 \gets \mathrm{compute\_increment}(u_0,\overline{u}_0,\{\mathcal{G}(u^{(n)})\}_n,\Gamma)$
\State $u_1 \gets \overline{u}_0 + \delta u_0$ 
\For{$j=1,\dots,J$}
\State $\lambda_j \gets \frac{j-1}{j+2},$
\State $v^{(n)}_j \gets u^{(n)}_j + \lambda_j (u^{(n)}_j -u^{(n)}_{j-1}), \forall n=1,\dots,N$
\State Compute $\overline{v}_j \gets \frac{1}{N}\sum_n v_j^{(n)}$
\State $\delta v_j \gets \mathrm{compute\_increment}(v_j,\overline{v}_j,\{\mathcal{G}(v^{(n)})\}_n,\Gamma)$
\State $u_{j+1} \gets \overline{v}_j + \delta v_j$ 
\EndFor
\State \Return $u_{J+1} = \{u^{(n)}_{J+1}\}_{n=1}^N$
\end{algorithmic}
\end{algorithm}

\begin{algorithm}
\caption{Nesterov-accelerated unscented Kalman inversion}
\label{alg:nest-uki}
\begin{algorithmic}[3]
\Function{$\mathrm{generate\_ensemble}$}{$m,C,r,\alpha,\Sigma_\omega$}
\State $\hat{m} \gets r + \alpha(m-r)$,  $\, \hat{C} \gets \alpha^2C + \Sigma_\omega$
\State $L \gets \mathrm{Cholesky}(\hat{C})$
\State $N \gets \mathrm{dim}(\hat{m})$, $\gamma\gets\sqrt{N}\min\{\sqrt{\frac{4}{N}},1\}$
\State $\begin{cases}
    u^{(0)} \gets \hat{m} & \\
    u^{(n)} \gets \hat{m} + \gamma \mathrm{col}_n(L)& n=1,\dots,N \\
    u^{(N+n)} \gets \hat{m} - \gamma \mathrm{col}_n(L)& n=1,\dots,N
\end{cases}$
\State \Return $u$
\EndFunction
\State 
\Function{$\mathrm{update\_mean\_cov}$}{$u,\{\mathcal{G}(u^{(n)})\}_n,y,\Sigma_\nu$}
\State Compute $\hat{m}$ and $\hat{C}$ from $u$,
\State $N \gets \mathrm{dim}(\hat{m})$ $\gamma\gets\sqrt{N}\min\{\sqrt{\frac{4}{N}},1\}$

\State $\hat{C}^{u\mathcal{G}} \gets \sum_n \frac{1}{2\gamma^2}(u^{(n)} - \hat{m})(\mathcal{G}(u^{(n)}) - \mathcal{G}(\hat{m}))^\top$
\State $\hat{C}^{\mathcal{G}\mathcal{G}} \gets \sum_n \frac{1}{2\gamma^2}(\mathcal{G}(u^{(n)}) - \mathcal{G}(\hat{m}))(\mathcal{G}(u^{(n)}) - \mathcal{G}(\hat{m}))^\top+\Sigma_\nu$ 
\State $m \gets \hat{m} + C^{u\mathcal{G}}(C^{\mathcal{G}\mathcal{G}})^{-1}(y-\mathcal{G}(\hat{m}))$ 
\State $C \gets \hat{C} - C^{u\mathcal{G}}(C^{\mathcal{G}\mathcal{G}})^{-1}(C^{u\mathcal{G}})^\top$
\State \Return{$m,C$}
\EndFunction
\State 
\Require $m_0, C_0, y, \mathcal{G}, r,\alpha,\Sigma_\omega, \Sigma_\nu$
\State $u_0\gets \mathrm{generate\_ensemble}(m_0,C_0,r,\alpha,\Sigma_\omega)$
\State $m_1,C_1 \gets\mathrm{update\_mean\_cov}(u_0,\{\mathcal{G}(u^{(n)}_0)\}_n,y, \Sigma_\nu)$

\For{$j=1,\dots,J$}
\State $u_j \gets \mathrm{generate\_ensemble}(m_j,C_j,r,\alpha,\Sigma_\omega)$
\State $\lambda_j \gets \frac{j-1}{j+2}$
\State $v^{(n)}_j \gets u^{(n)}_j + \lambda_j(u^{(n)}_j - u^{(n)}_{j-1})$

\State $m_{j+1},C_{j+1} \gets\mathrm{update\_mean\_cov}(v_j,\{\mathcal{G}(v^{(n)}_j)\}_n,y,\Sigma_\nu)$
\EndFor
\State \Return{$m_{J+1},C_{J+1}$}
\end{algorithmic}
\end{algorithm}

\section{Additional experiments}\label{sec:additional_experiments}

We briefly present three additional experiments in support of the above results.

\begin{figure}
    \centering
    \includegraphics[width=0.99\linewidth]{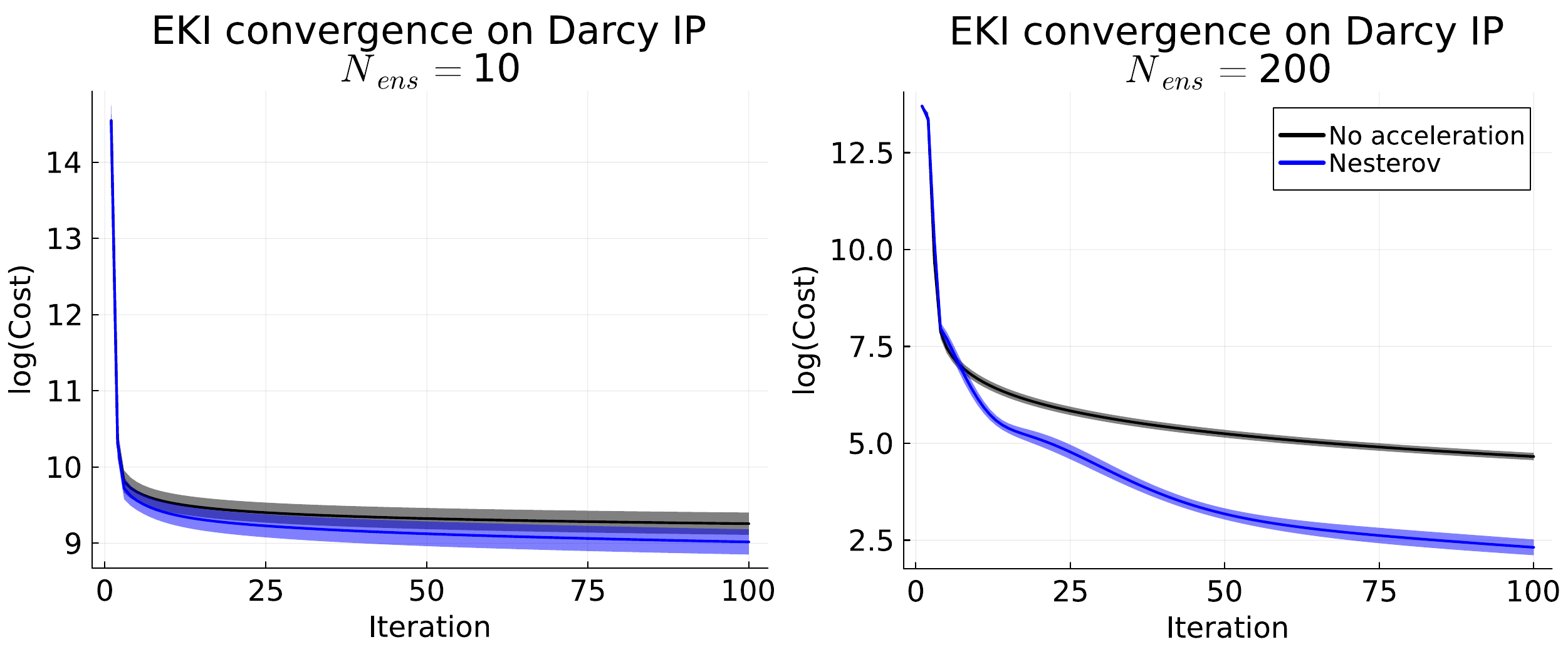}
    \caption{The effect of ensemble size on EKI convergence with and without acceleration, on the Darcy inverse problem of Section\ref{sec:ex-ip}, as measured by the logarithm of the cost function \eqref{eq:cost}. Nesterov-accelerated EKI applied to the Darcy problem, with ensemble sizes 10 and 200. An experiment with ensemble size 52 is given in Figure \ref{fig:nesterov_results}.}
    \label{fig:ens_sizes}
\end{figure}
We test the performance of Nesterov-accelerated EKI on the Darcy problem with two different ensemble sizes, $N=10$ and $N=200$, shown in Figure \ref{fig:ens_sizes}. In the case of $N=10$, both algorithms can explore only the nine-dimensional subspace spanned by the initial ensemble, while the true parameters are in a fifty-dimensional space; this leads to higher costs and large spread across trials due to sampling error. For $N=200$ case there is no restriction as the ensemble subspace contains the true parameters, and the cost is greatly reduced. In both cases, we again observe clearly that the Nesterov acceleration is beneficial for convergence, though the spread due to sampling error in the small-ensemble regime causes some overlapping performance in the first few iterations.



\begin{figure}
    \centering
    \includegraphics[width=0.9\linewidth]{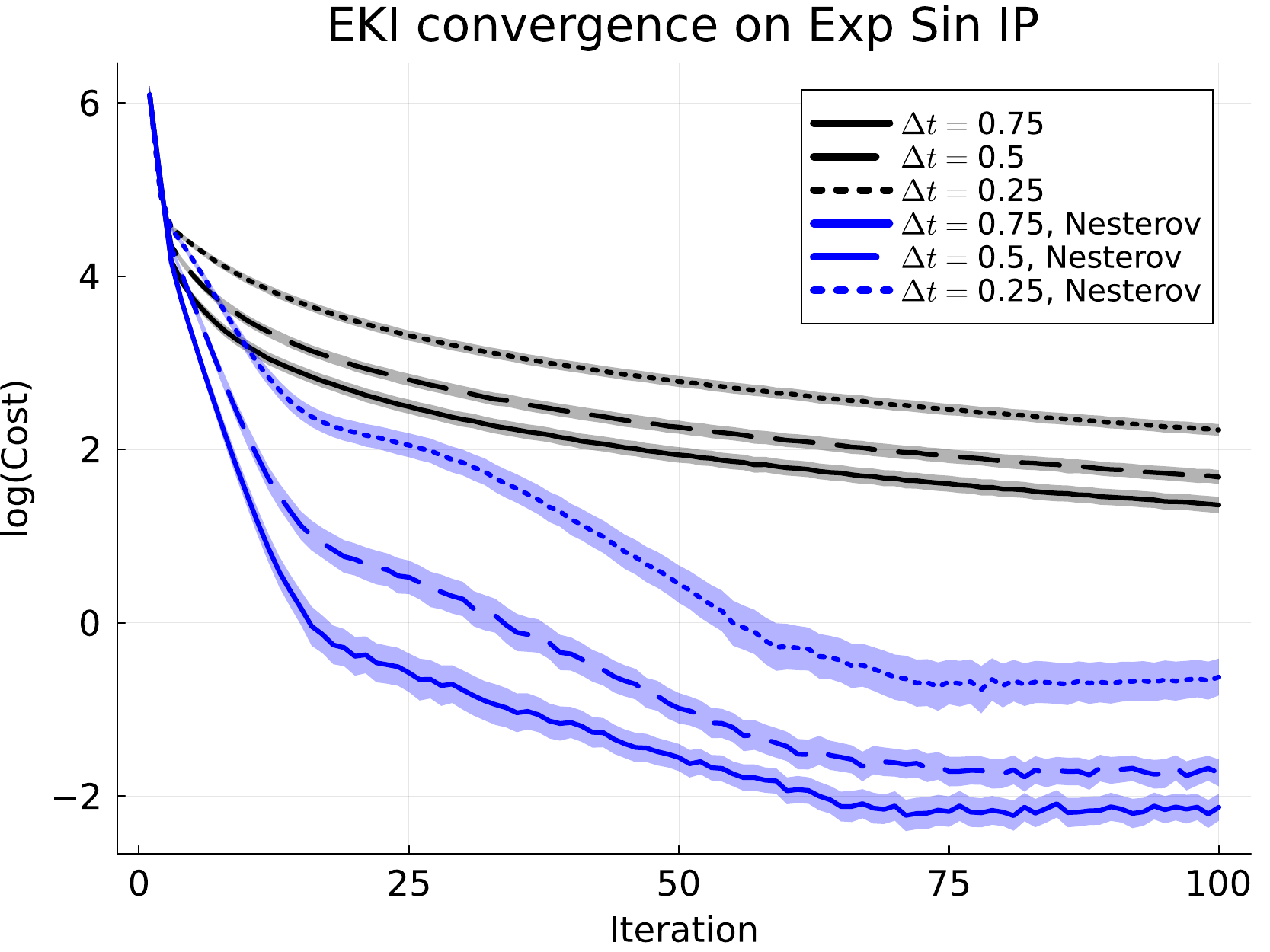}
    \caption{The effect of taking different $\Delta t$ on EKI convergence with and without acceleration on the Exp Sin inverse problem of Section \ref{sec:ex-ip}, as measured by the logarithm of the cost function \eqref{eq:cost}.}
    \label{fig:nesterov_timesteps}
\end{figure}
Figure~\ref{fig:nesterov_timesteps} shows the results of applying EKI and Nesterov-accelerated EKI on the exponential sine problem with varying $\Delta t$. Although the timestep has an influence on the performance, the Nesterov-accelerated EKI with the worst choice of timestep still outperforms (after 10 iterations) the non-accelerated EKI with the best choice of timestep. This suggests that the Nesterov acceleration can compensate for a poorly chosen timestep.

\bibliographystyle{elsarticle-num} 
\bibliography{references.bib}



\end{document}